# A new criterion based on Kullback-Leibler information for space filling designs

Astrid Jourdan[a,*], Jessica Franco[b]

[a] Department of mathematics, E.I.S.T.I, 26 avenue des Lilas, 64062 Pau cedex 9, France

[b] Total- DGEP/GSR/TG/G&I, avenue Larribau, 64018 Pau Cedex – France

---

ABSTRACT

Experimental designs are tools which can drastically reduce the number of simulations required by time-consuming computer codes. One strategy for selecting the values of the inputs, whose response is to be observed, is to choose these values so that they are spread evenly throughout the experimental region, according to "space filling designs". In this article, we suggest a new criterion based on the Kullback-Leibler information for design construction. The aim is to minimize the difference between the empirical distribution of the design points and the uniform distribution which is equivalent to maximizing the Shannon entropy. The entropy is estimated by a Monte Carlo method, where the density function is replaced with its kernel density estimator or by using the nearest neighbor distances.

---

## 1. Introduction

For many scientific phenomena, physical experimentation is very expensive, time-consuming, or impossible. Engineers and scientists have been in the forefront of developing mathematical models and numerical solutions to describe physical systems. As models become more sophisticated, running times increase rapidly, and computer experiments are necessary to characterize physical phenomena.

* Corresponding author. Tel. +33 (0)5 59 14 85 33 ; Fax. +33 (0)5 59 14 85 31

*E-mail addresses* : astrid.jourdan@eisti.fr (A. Jourdan), jessica.franco@total.com (J. Franco)

A computer experiment consists of running a simulation with an input vector x which specifies the values of the governing input parameters (input factors) of the computer model. The outputs, $y(x_1),...,y(x_n)$, at a given set of inputs $x_1,...,x_n$, are used to provide a predictor of the simulated response. In this paper, we suppose that no information is available about the relationship between the response and the input factors (exploratory phase). Hence, the designs should allow one to adapt a variety of statistical models (kriging models, neural network models, ...). One strategy consists of selecting the set of inputs in order to fill up the experimental region in a uniform fashion. These kinds of designs are called “space filling designs” and have been widely investigated in the past decade (e.g. Koehler & Owen (1996) and Franco (2008)). They provide information about all parts of the experimental region and enable one to spot possible irregularities of the computer response within the experimental region.

In this paper, we propose a new criterion based on the Kullback-Leibler information (KL information) for design construction. As with the discrepancy method, the KL information measures the difference between the empirical distribution of the design points and the uniform distribution. The idea is to minimize this difference by using an exchange algorithm.

In section 2, the properties of the KL information are presented, especially the link with the Shannon entropy in the case of the uniform distribution. In sections 3 and 4, two methods of the entropy estimation are investigated. The first one is a Monte Carlo method where the density function is replaced with its Gaussian kernel estimate. The second one avoids the density estimation by using the nearest neighbor distances. In section 4, we compare the new designs with the most commonly used space filling designs.

## 2. Kullback-Leibler information

Suppose that the design points $X_1,...,X_n$, are $n$ independent observations of the random vector $X=(X^1,...,X^d)$ with absolutely continuous density function $f$ concentrated on the unit cube $[0,1]^d$. The aim is to select the design points in such way as to have the density function “close” to the uniform density function. The Kullback-Leibler (KL) information measures the difference between two density functions $f$ and $g$ supported by E,

$$I(f,g) = \int_E f(x)\ln\left(\frac{f(x)}{g(x)}\right)dx\ .$$

Note that the KL information is not a metric function since the triangle inequality is not satisfied and the information is not symmetrical. Hence, the density functions can not be exchanged and *f* is the “design” density and *g* is the “target” density, *i.e* the uniform density. However, the KL information is always non-negative, $\mathrm{I}(f,g) \geq 0$, with $\mathrm{I}(f,g) = 0$ if and only if *f*=*g*. The aim is then to minimize the KL information so as to make *f* converge towards the uniform density.

If *g* is the uniform density function, then the KL information becomes

$$\mathrm{I}(f) = \int f(\mathrm{x})\ln\left(f(\mathrm{x})\right)\mathrm{dx} = -\mathrm{H}\left[f\right],$$

where H[*f*] denotes the Shannon entropy. Therefore, minimizing the KL information amounts to maximizing the entropy.

We recognize the definition of the maximum entropy sampling commonly used in computer experiments (Shewry & Wynn (1987), Currin *et al.* (1988)). The difference lies in the fact that the entropy used in this paper is independent of any metamodel. The aim of maximizing the entropy is not the gain in information on the model parameters. On the other hand, for *f* concentrated on $[0,1]^{d}$, one always has H(*f*)≤0 and the maximum value of H(*f*), zero, being uniquely attained by the uniform density. This latter property confirms our choice to minimize the KL information. Note that Dudewicz *et al.* (1981) proposed an entropy-based test of uniformity but only for *d*=1. In order to avoid any confusion with the maximum entropy designs mentioned above, the designs constructed in this paper are called «minimal KL information designs ».

Finally, we use an exchange algorithm that maximizes the entropy for the design construction. The main point of the method is the estimation of the Shannon entropy. Beirlant *et al.* (1997) gave an overview of several techniques of estimation. Two methods have been selected. The first one is a Monte Carlo method where the density function is replaced with its Gaussian kernel estimate. The second one is based on the nearest neighbor distances.

## 3. Estimate of entropy based on the Monte Carlo method (MC)

Given that the entropy can be written as follows,

$$\mathrm{H}(\mathrm{X}) = -\mathrm{E}\left[\ln\left(f(\mathrm{X})\right)\right]$$

where E denotes the expectation for *f*, the Monte Carlo method provides an unbiased and consistent estimate of the entropy,

$$\hat{H}(X) = -\frac{1}{n}\sum_{i=1}^{n} \ln f(X_i) \,.$$

Ahmad and Lin (1976) proposed to replace the unknown density function *f* with its kernel density estimate (Silverman 1986, Scott 1992),

$$\hat{f}(x) = \frac{1}{nh^d}\sum_{i=1}^{n} K\left(\frac{x - X_i}{h}\right),\ \forall x \in [0,1]^d$$

The quality of a kernel estimate depends essentially on the value of its bandwidth *h* (smoothing parameter). In our application, the bandwidth is chosen using Scott's rule where the standard deviation estimates are replaced with the standard deviation of the uniform distribution on [0,1],

$$\hat{h} = \frac{1}{\sqrt{12}}\frac{1}{n^{1/(d+4)}}$$

The choice of kernel function K is much less important for the behavior of the estimate than the choice of h. Most of the common kernels (uniform, Epancehnikov, triangle, ...) have a bounded support (unit sphere), so that, in our application, the probability that the kernel values are not zero is extremely small. So, the chosen kernel function K is the multivariate Gaussian distribution,

$$K(z) = \frac{(2\pi)^{-d/2}}{s^d}\exp\left[-\frac{1}{2s^2}\|z\|^2\right],$$

where $s^2$ is the variance of the uniform distribution on [0,1] multiplied by the dimension d, $s^2 = \frac{d}{12}$.

Joe (1989) obtained asymptotic bias and variance terms for the estimator described above. The bias depends naturally on the size *n* and the dimension *d* but also on the bandwidth *h* (fixed during the exchange algorithm). Moreover, he pointed out that the sample size needed for good estimates increases rapidly with the dimension (fig 1a). Hence, in the context of experimental design, the entropy estimation will be approximate since the size *n* is small. However, the objective is not to compute an accurate approximation of the entropy. And, despite the low accuracy, the exchange algorithm converges rapidly (fig 1b) toward designs with the expected properties (filling up the space) (fig. 3).

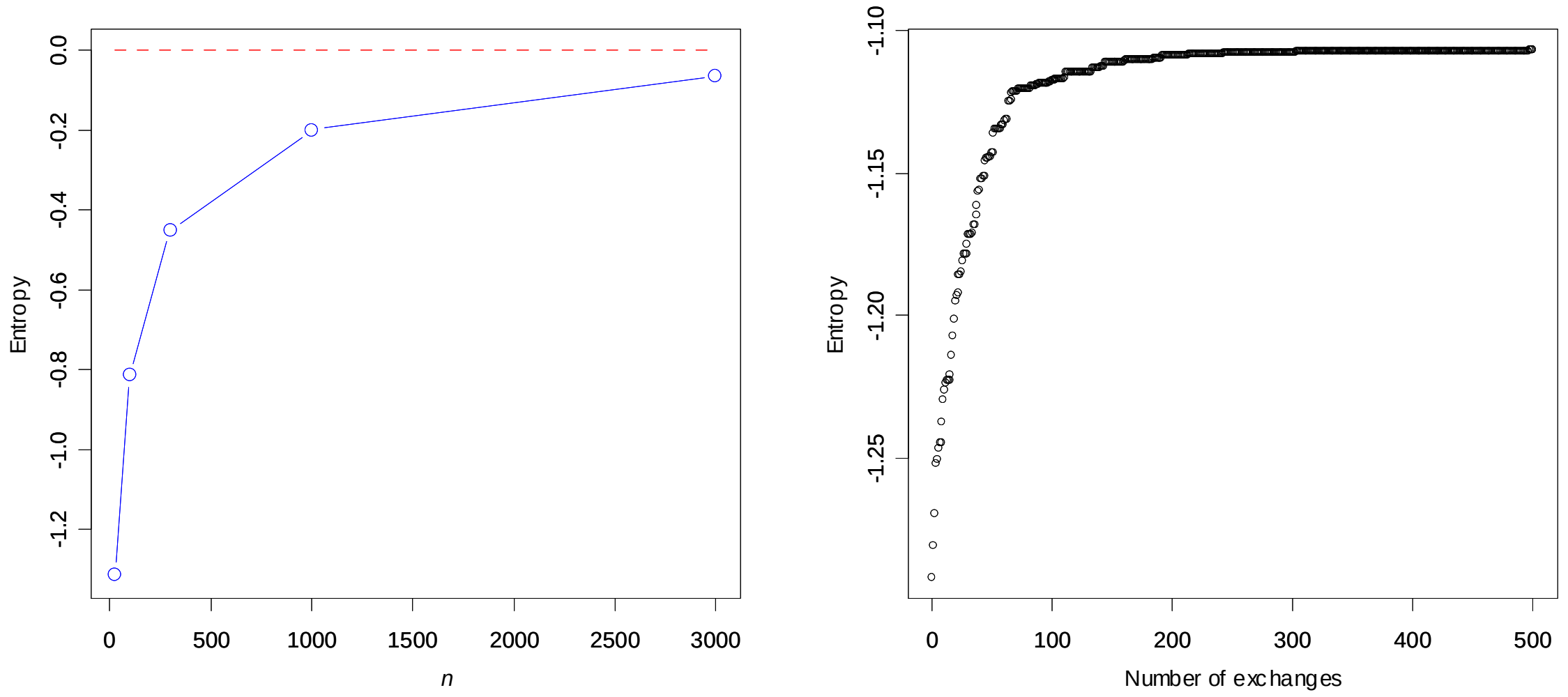

(a) Entropy (MC) *vs.* the size *n*

(b) Entropy (MC) *vs.* the number of exchanges (*n*=30)

**Fig. 1.** Convergence of the Monte Carlo method (*d*=3)

## 4. Estimate of entropy based on nearest neighbor distances (NN)

The performance of the Monte Carlo methods generally reduces as dimensionality increases. Moreover, the density function needs to be estimated for each exchange in the algorithm. Hence, the design construction with the previous method becomes time-consuming as the sample size increases. Kozachenko and Leonenko (1987) proposed estimating the entropy by

$$\hat{H}(f) = \frac{d}{n}\sum_{i=1}^{n} \ln \rho_i + \ln\left[\frac{\pi^{d/2}}{\Gamma(d/2+1)}\right] + C_E + \ln(n-1)$$

where $C_E \approx 0.5772$ is the Euler constant, $\Gamma$ is the Gamma function and $\rho_i$ is the nearest neighbor distance of $X_i$ and the other $X_j$,

$$\rho_i = \min_{j \neq i, 1 \leq j \leq n} \left\| X_i - X_j \right\|.$$

They proved the asymptotic unbiasedness and the consistency of the estimator. Figure 2a illustrates the consistency and the bias since the estimated values are positive whereas the entropy is negative. Figure 2b shows that the exchange algorithm converges more slowly than the previous one. Finally, even if the density function estimation is avoided, this method of construction requires nearly the same running-time.

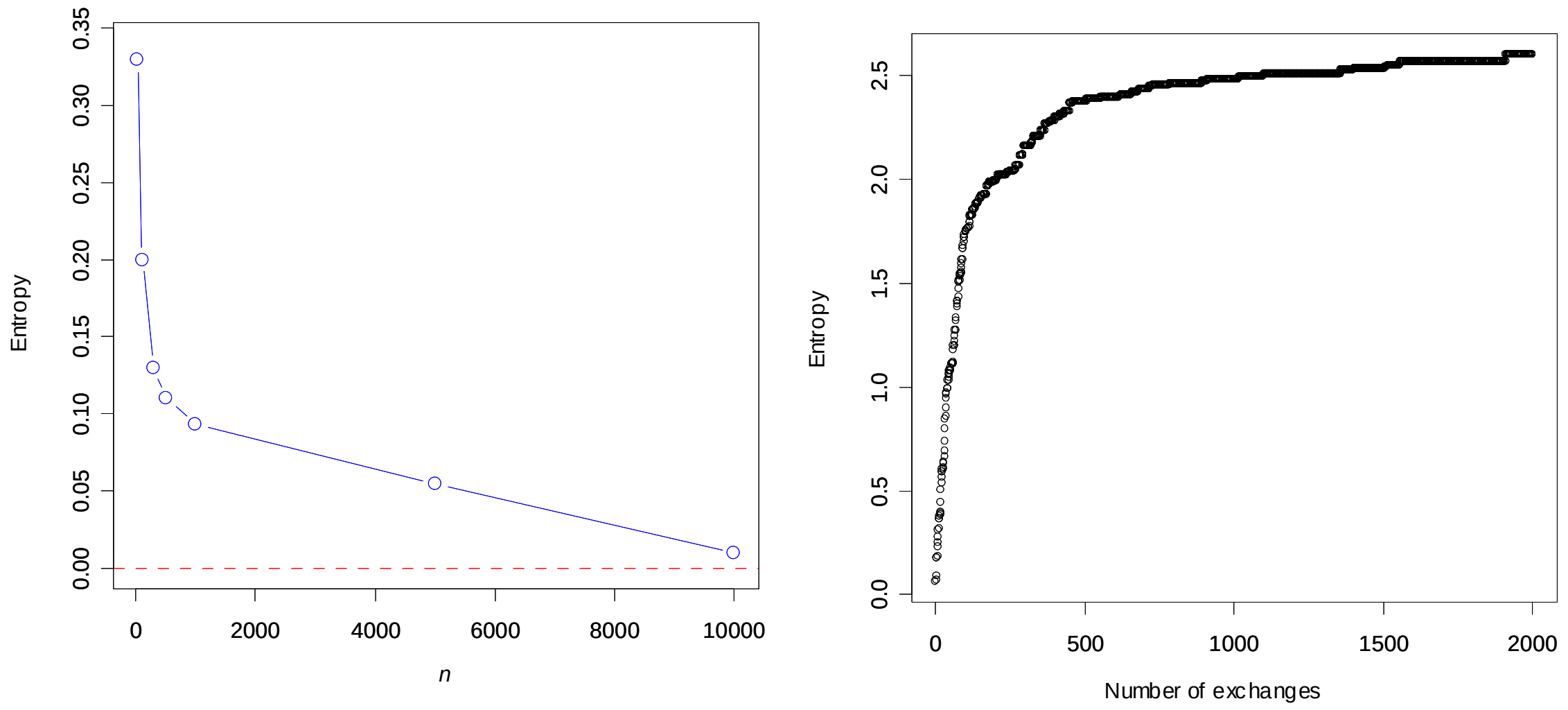

(a) Entropy (NN) *vs*. the size *n*

(b) Entropy (NN) *vs*. the number of exchanges (*n*=30)

**Fig. 2.** Convergence of the nearest neighbor method (*d*=3)

Remark : Singh *et al.* (2003) generalized the method to the $k^{th}$ nearest neighbor in order to reduce the high fluctuations of the estimation due to the small fluctuations of the small $\rho_i$ values. This method is justified in the context of real experiments where the $\rho_i$ values are not recorded to high accuracy. Conversely, the values are exact in simulation and it is preferable to emphasize the small distances.

## 5. Design comparison

The aim of this section is to study the properties of the designs constructed by the methods described above, and to compare them with the usual space filling designs.

In an exchange algorithm, the design depends more or less on the initial setting. Consequently, several initializations are run, and the best design is selected so as to reduce the dependence. However, in this section, several designs are built from only one initialization in order to study this dependence.

### *5.1. Improvement of the initial setting*

The aim of the construction methods is to fill up the parameter space in a uniform fashion. It can be seen on figure 6 that the objective is reached for the two methods. Whatever the initialization, the algorithm converges toward designs with the same characteristics :

- The points lie on the edge of the factor space but also in the interior like a scrambled regular grid (quasi-periodical distribution). Such distribution assures the points are spread evenly in the unit cube.
- Contrary to a regular grid, many levels are tested for each parameter.
- For small size n, the minimal KL designs will generally lie on the exterior of $[0,1]^d$ and fill in the interior as n becomes large. Koelher & Owen (1996) pointed out the same behavior for maximin designs.

The MST criterion (§ 5.3.) shows that the designs keep these characteristics in dimension d≥2.

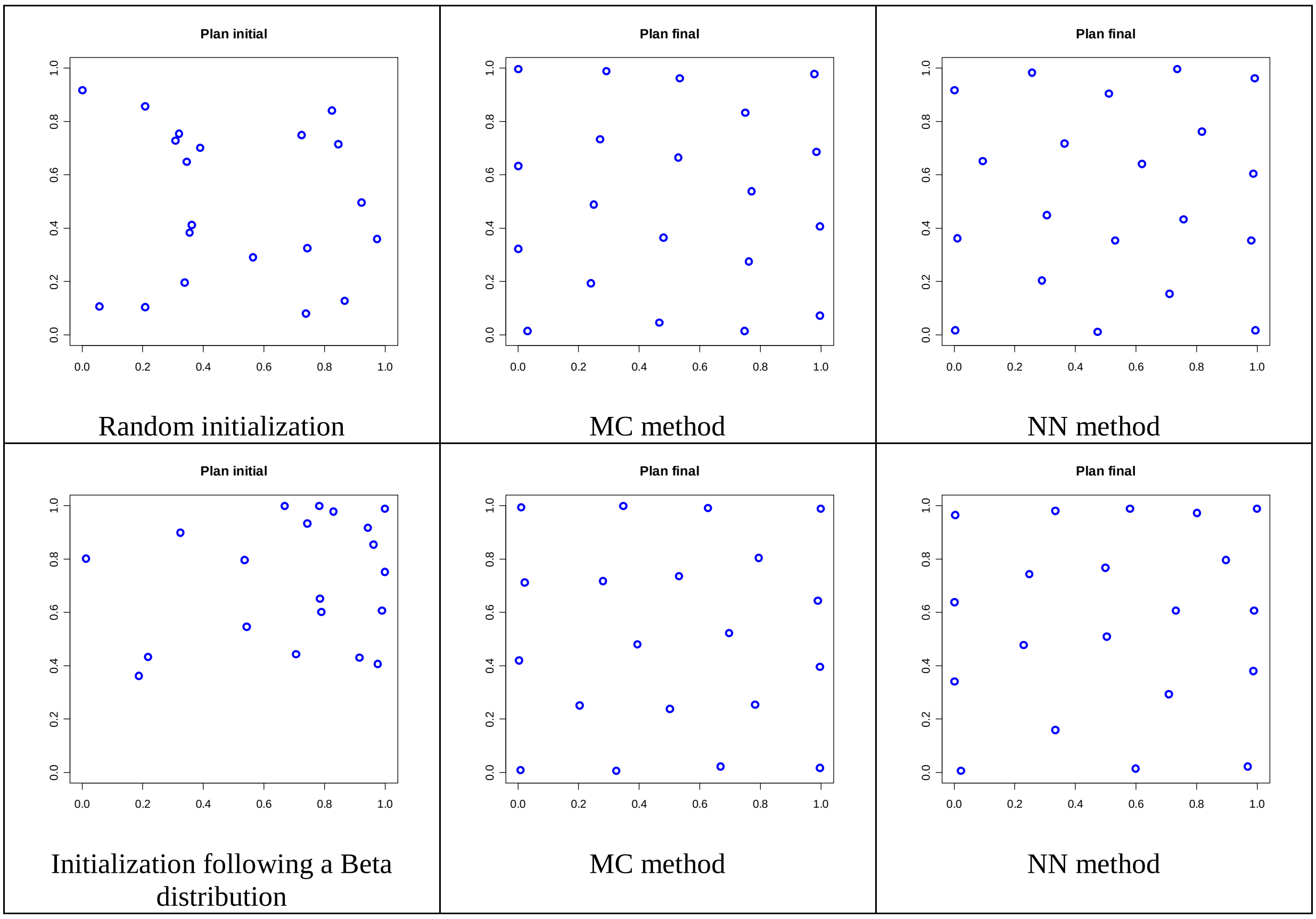

**Fig. 3.** Results of the exchange algorithm for designs of size 20 in dimension 2

### 5.2. *Comparison with traditional space filling designs*

**Criteria**

It is impossible to study visually the point distribution when d≥2. The analysis requires criteria to assess the uniformity and/or the distribution throughout the unit cube. The criteria selected in this paper are intrinsic, that is, they are independent of any kind of metamodel.

- The cover measure (Cov) calculates the difference between the design and a regular grid (Gunzburger *et al.* (2004)). The goal is to minimize the criterion in order to be close to a regular grid and, in so doing, to guarantee that the points fill up the space.
- Johnson *et al.* (1990) proposed the maximin and minimax distances to spread the points evenly in the unit cube. The maximin criterion (Mindist) maximizes the minimal distance between the points.
- Contrary to the previous criterion, the discrepancy is not based on the distance between the design points. This criterion measures the difference between the empirical cumulative distribution and the uniform cumulative distribution (Niederreiter (1987)). Two kinds of discrepancy are selected, the L2-discrepancy (DL2) and the centered L2-discrepancy (DC2).

**Traditional space filling designs**

The designs selected for the comparison are the most common designs used for computer experiments. Koehler & Owen (1996) and Franco (2008) proposed an overview of space filling designs.

- Random : designs built by simple applications of the random function.
- LH : latin hypercube designs (without optimization).
- Discr : sequences with small discrepancy (Sobol, Niederreiter, Hammersley, Halton) (Niederreiter (1987)).
- Dmax : designs which maximize the determinant of a covariance matrix (Shewry & Wynn (1987), Currin *et al.* (1988)). These designs are commonly used when the metamodel is a kriging surface.
- Strauss : designs created using a Strauss procedure, which considers the repulsion between two points, in such a way as to maximize filling the space of the parameters (Franco, 2008).
- Maximin : optimal designs based on the maximin distance.
- MC : minimal KL information designs constructed by the Monte Carlo method.
- NN : minimal KL information designs constructed by the nearest neighbor method.

The minimal KL designs provide the best results except for the discrepancy (fig. 4 and fig. 5). This conclusion is surprising since the two criteria, discrepancy and KL information, minimize the difference between the point distribution and the uniform distribution, whereas the cover measure and the maximin criteria are based on the distance between the design points. This result may be explained by the difficulty in evaluating the discrepancy.

Figure 5 shows that the design qualities do not decrease when the dimension $d$ increases. The minimal KL information designs compete with maximin designs (even concerning the maximin criteria) which are commonly used in the exploratory phase.

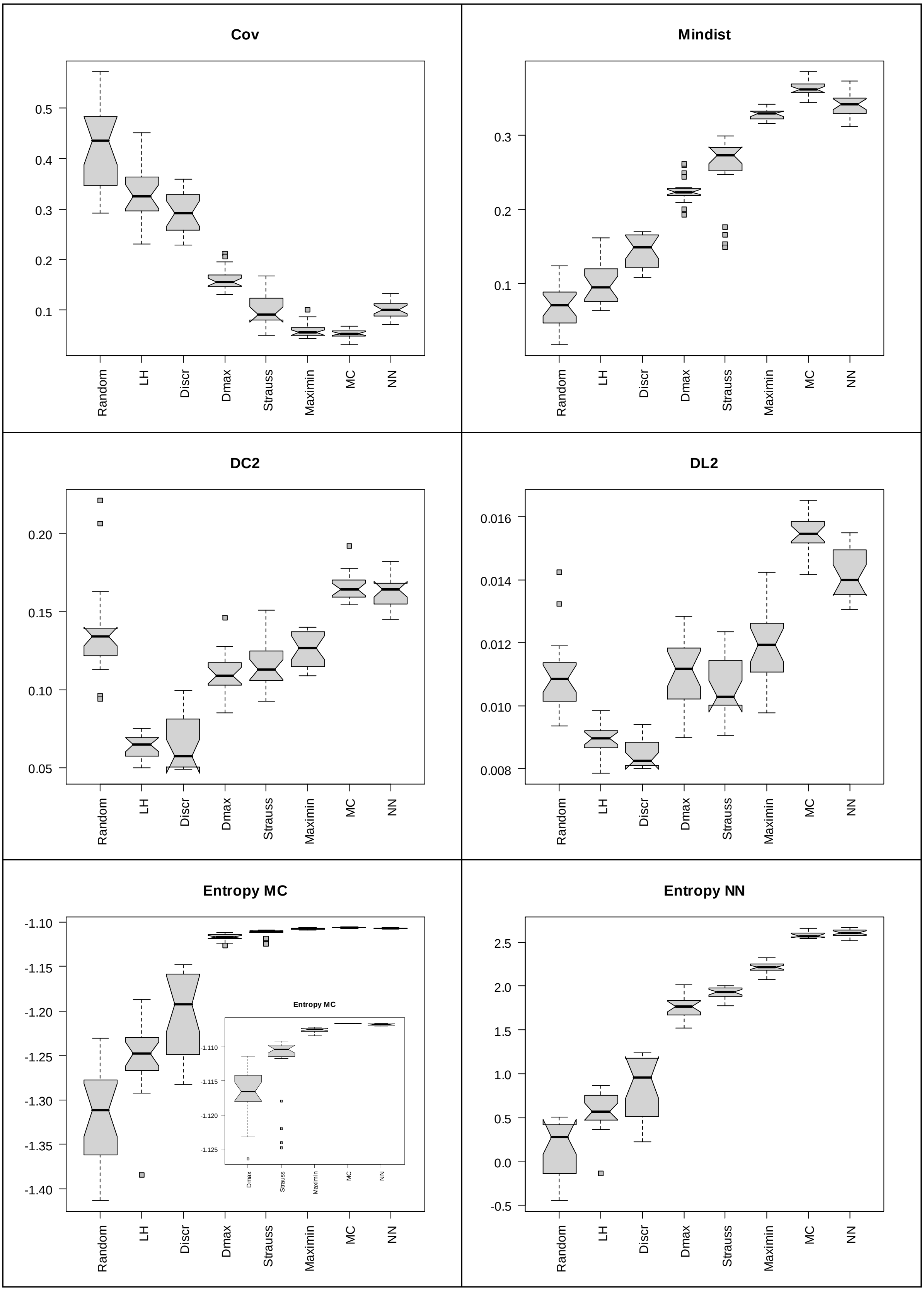

**Fig. 4.** Criterion for 20 designs of size 30 with dimension 3

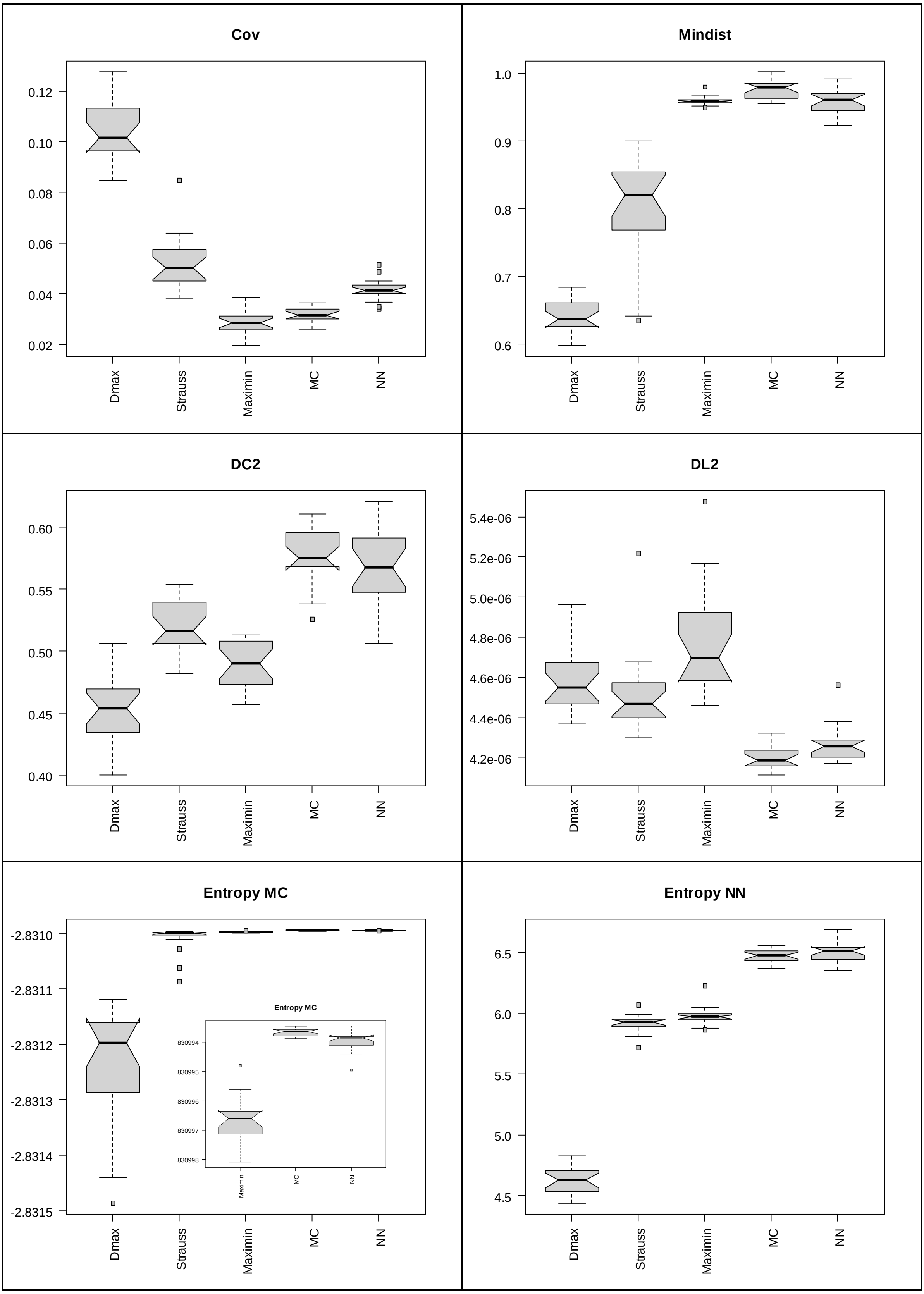

**Fig. 5.** Criterion for 20 designs of size 100 with dimension 10

*5.3. MST criterion*

The MST (Minimal Spanning Tree) criterion proposed by Franco *et al.* (2008) (following the studies of Wallet *et al.*, 1998) classifies the designs according to their structure: grid, cluster,… . The classification depends on the mean and standard deviation of the edge lengths of the MST associated with the design under study (fig. 6). The designs of the region of quasi-periodicity represent the best compromise between regular grid (space filling) and uncertain distribution (uniformity). The designs of minimal KL information also give the best results in this case (figures 7 and 8) with a clear demarcation in dimension 3. This result confirms the remarks on the characteristics of the designs of dimension 2 (§5.1). The designs constructed using the Monte Carlo method with Gaussian kernels are still higher performing than those constructed by nearest neighbor method.

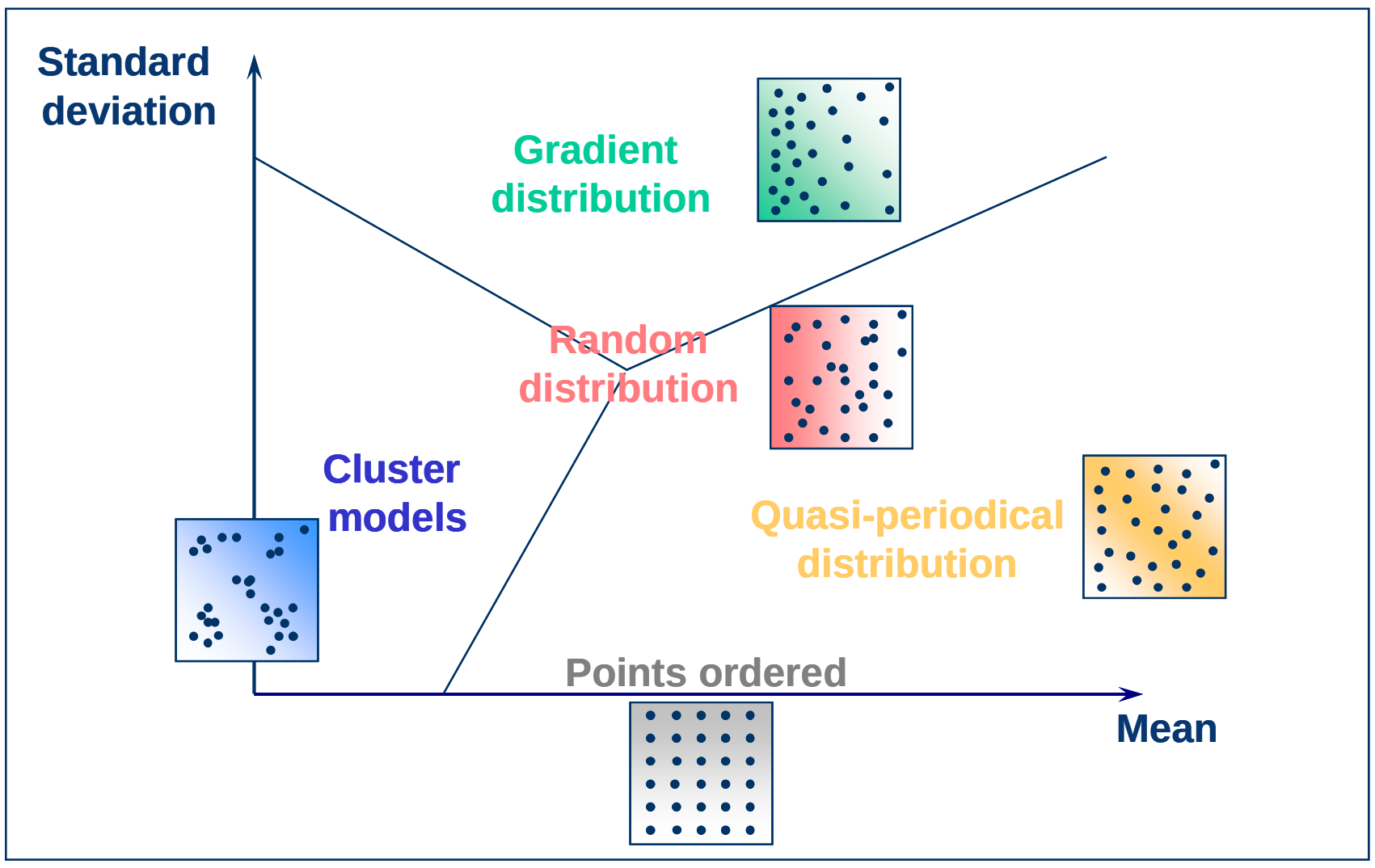

**Fig. 6.** Empirical distribution created using MST criteria

The classification depends on the mean and standard deviation of the edge lengths of the MST associated with the design under study. The designs of the region of quasi-periodicity represent the best compromise between regular grid (space filling) and uncertain distribution (uniformity). The designs of minimal KL information also give the best results in this case (fig. 7 and fig. 8) with a clear demarcation in dimension 3. This result confirms the remarks on the characteristics of the designs of dimension 2 (§5.1). The designs constructed using the Monte Carlo method with Gaussian kernels are still higher performing than those constructed by nearest neighbor method.

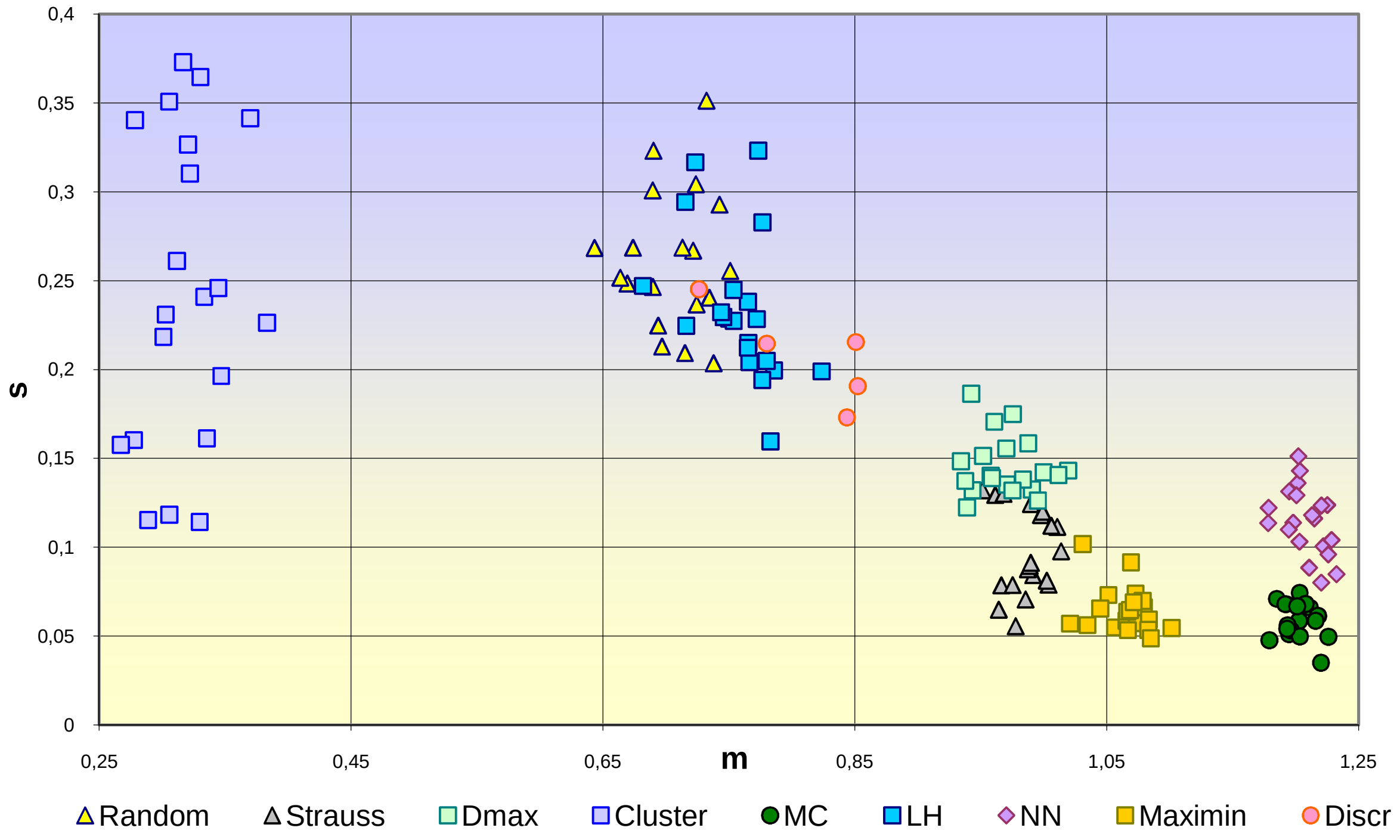

**Fig. 7.** Representation of the MST criterion calculated on 20 designs with 30 points in dimension 3

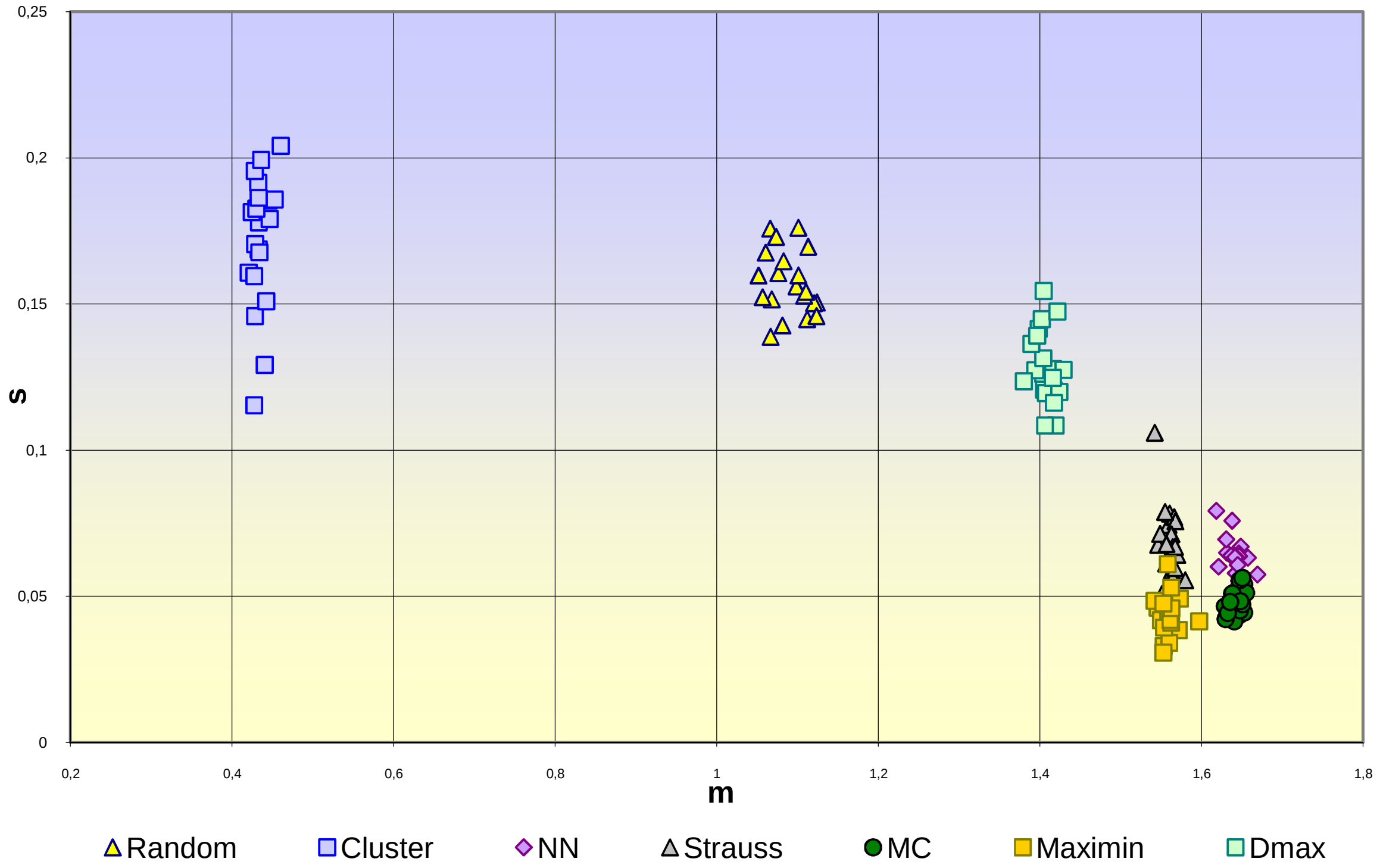

**Fig. 8.** Representation of the MST criterion calculated on 20 designs with 100 points in dimension 10

## 6. Conclusion

The aim of the criteria presented in this paper is to minimize the difference between the distribution of the design points and the uniform one by using the KL information. The criterion is equivalent to maximizing the Shannon entropy, and the design construction is then based on the entropy estimation. Two methods have been investigated, and they provide roughly the same results. The first one is a Monte Carlo method where the density function is replaced with its Gaussian kernel estimation. The second one avoids the density estimation by using the nearest neighbor distances. The tests show that minimal KL information designs achieve the main goal, which is to spread the design points evenly throughout the experimental region. They are indisputably the best designs with regard to the usual criteria, even in high dimensions. They compete with maximin designs which are widely used in the exploratory phase.

Despite the context of experimental designs, which is far-removed from asymptotic conditions, a criterion based on the Shannon entropy estimation provides very good results. However, the Monte Carlo method performs less well as the dimension increases and the design size decreases. Moreover, the Monte Carlo method requires running time, since the density has to be estimated for each exchange. Leonenko *et al.* pointed out that nearest-neighbor methods show promise when the dimension is large, but the tests show a slow convergence of the exchange algorithm. One idea consists of replacing the Shannon entropy by another one, for instance, the Rényi entropy (estimated by minimal spanning trees, Hero *et al.*, 2002), or the Tsallis entropy (analytic expression when the density is replaced with its kernel estimate, Bettinger *et al.*, 2008).

**Acknowledgment**

This work was supported by Total. We are grateful for the encouragements from Bernard Corre.

**References**

Ahmad I.A., Lin P.E., 1976. A nonparametric estimation of the entropy for absolutely continous distributions. IEEE Trans. Information Theory, 35, pp. 688-692.

Ahmed N.A, Gokhale D.V., 1989. Entropy expressions and their estimators for multivariate distributions. IEEE Transactions on Information Theory, 35(3), pp. 688-692.

Beirlant J., Dudewicz E.J., Györfi L., Van Der Meulen E.C., 1997. Nonparametric entropy estimation : an overview. Int. J. Math. Stat. Sci., 6(1) 17-39.

Bettinger R., Duchêne P., Pronzato L., Thierry E., 2008. Design of experiments for response diversity. In Proc. 6th International Conference on Inverse Problems in Engineering (ICIPE), Journal of Physics: Conference Series}, Dourdan (Paris), 15-19 juin 2008, to appear. http://hal.archives-ouvertes.fr/hal-00290418/fr/

Chaloner K., Verdinelli I., 1995. Bayesian Experimental Design: A review. Statist. Sci., 10, 237-304.

Currin C., Mitchell T., Morris M., Ylvisaker D., 1988. A bayesian approach to the design and analysis of computer experiments. ORNL Technical Report 6498, available from the national technical information service, Springfield, Va. 22161.

Dimitriev Y.G, Tarasenko F.P., 1973. On the estimation functions of the probability density and its derivatives. Theory Probab. Appl., 18, 628-633.

Dudewicz E.J., Van Der Meulen E.C., 1981. Entropy-Based tests of uniformity. J. Amer. Statist. Assoc., 76, 967-974.

Franco J., 2008. Planification d'expériences numériques en phase exploratoire pour des codes de calculs simulant des phénomènes complexes. Thèse présentée à l'Ecole Nationale Supérieure des Mines de Saint-Etienne

Franco J., Vasseur O., Corre B., Sergent M., 2008. Minimum Spanning Tree : A new approach to assess the quality of the design of computer experiments. To appear in Chemometrics and Intelligent Laboratory Systems.

Gunzburger M., Burkardt J., 2004. Uniformity measures for point sample in hypercubes. https://people.scs.fsu.edu/~burkardt/pdf/ptmeas.pdf

Harner, E. James, Singh H., Li, S., and Jun Tan, 2004. Computational Challenges in Computing Nearest Neighbor Estimates of Entropy for Large Molecules. Proceeding of the 36th Symposium on the Interface: Computational Biology and Bioinformatics, may 26-29.

Hero A., Bing Ma, Michel O., Gorman J., 2002. Applications of entropic spanning graphs. Signal Processing Magazine, IEEE, 19, 85 – 95.

Joe H., 1989. Estimation of entropy and other functional of multivariate density. Ann. Int. Statist. Math., 41, 683-697.

Johnson M.E., Moore L.M., Ylvisaker D. (1990). Minimax and maximin distance design. J. Statist. Plann. Inf., 26,131-148.

Koehler J.R., Owen A.B, 1996. Computer Experiments. Handbook of statistics, 13, 261-308.

Kosachenko L.F., Leonenko N.N., 1987. Sample estimate of entropy of a random vector. Problem of Information Transmission, 23, 95-101.

Kullback S., Leibler R.A., 1951. On information and sufficiency. Ann. Math. Statist., 22 79-86.

Leonenko N, Pronzato L, Savani V. A class of Rényi information estimators for multidimensional densities. Annals of Statistics , to appear.

Niederreiter H., 1987. Point sets and sequences with small discrepancy. Monasth. Math., 104, 273-337.

Scott D.W., 1992. Multivariate Density Estimation : Theory, practice and visualization, John Wiley & Sons, New York, Chichester.

Sebastiani P. & Wynn H.P., 2000. Maximum entropy sampling and optimal Bayesian experimental design. J. Royal Statist. Soc., 62, 145-157

Shewry M.C., Wynn H.P., 1987. Maximum Entropy Sampling. J. Appl. Statist., 14, 165-170.

Silverman B.W., 1986. Density estimation for statistics and data analysis. Chapman & Hall, London.

Singh, H., Misra, N., Hnizdo, V., Fedorowicz, A. & Demchuk, E., 2003. Nearest neighbor estimates of entropy. Am. J. Math. Manage. Sci., 23, pp. 301-321.

Stein M. (1987). Large sample properties of simulations using latin hypercube sampling. Techometrics, 29, 143-151.

Thiémard E. (2000). Sur le calcul et la majoration de la discrépance à l’origine. Thèse présentée au département de mathématiques de l’école polytechnique fédérale de Lausanne

Wallet F., Dussert C. (1998). Comparison of spatial point patterns and processes characterization methods. Europhysics Lett., 42, 493-498.